# The Voronoi Cell in a saturated Circle Packing

### and an elementary proof of Thue´s theorem

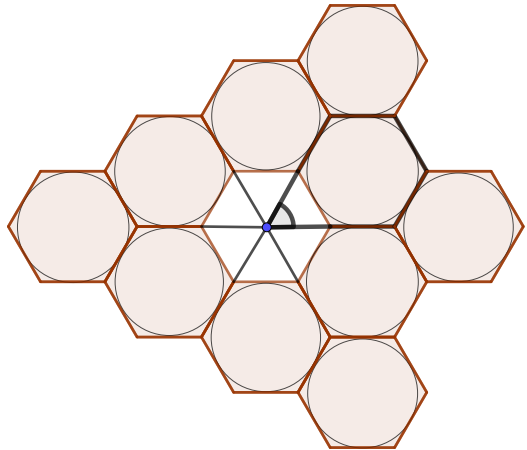

*Voronoi Cell in a non-saturated hexagonal packing*

Dr. Max Leppmeier
  Lehrstuhl für Mathematik und ihre Didaktik
  Universität Bayreuth
  max.leppmeier@uni-bayreuth.de

Max Leppmeier studied mathematics and physics at the Ludwig Maximilian University of Munich. He wrote the book "Kugelpackungen von Kepler bis heute" ("Sphere Packings from Kepler until today") and received his doctorate (Ph.D.) in mathematics and its didactics from the University of Bayreuth. Recently his book "Mathematische Begabungsförderung am Gymnasium" ("Mathematical Talent Development – Concepts for Teaching and School Development") has been published.

**Abstract**
The famous Kepler conjecture has a less spectacular, two-dimensional equivalent: The theorem of Thue states that the densest circle packing in the Euclidean plane has a hexagonal structure. A common proof uses Voronoi cells and analyzes their area applying Jensen´s inequality on convex functions to receive a local estimate which is globally valid.
Based on the concept of Voronoi cells, we will introduce a new tessellation into so-called L-triangles which can be related to fundamental parallelograms of lattice circle packings. Therefore a globally disordered circle packing can be reduced to locally ordered configurations: We will show how the theorem of Lagrange on lattice circle packings can be applied to non-lattice circle packings. Thus we receive a new proof of Thue´s theorem.

**Zusammenfassung**
Die berühmte Kepler-Vermutung hat ein Pendant im Zweidimensionalen: Der Satz von Thue besagt, dass die dichteste Kreispackung in der Euklidischen Ebene eine hexagonale Struktur aufweist. Ein Standardbeweis analysiert die Fläche von Voronoi-Zellen mit Hilfe der Jensen´schen Ungleichung für konvexe Funktionen. Dies ergibt eine lokale Abschätzung, die global gültig ist.
Aufbauend auf dem Begriff der Voronoi-Zelle führen wir eine neue Zerlegung in sog. L-Dreiecke ein, die in Bezug zu Fundamentalparallelogrammen von Kreisgitterpackungen gesetzt werden können. Auf diese Weise kann eine global ungeordnete Kreispackung zurückgeführt werden auf eine lokal geordnete Konfiguration: Es wird gezeigt, wie der Satz von Lagrange über dichteste Kreisgitterpackungen auch auf nicht gitterförmige Kreispackungen angewandt werden kann. Auf diese Weise erhalten wir einen neuen Beweis für den Satz von Thue.

# 1 Introduction

The Kepler conjecture was introduced by Kepler in 1611 and asserts that a densest packing of equal spheres in three-dimensional Euclidean space is given by the „cannonball" packing, or face-centered-cubic (fcc) packing [1, p. 5] [2] [3] [4] [5] [6]. The fcc packing can be described either by its lattice and fundamental parallelotope or by its reciprocal lattice and the Voronoi cell. Hales is regarded as the conqueror of the Kepler problem [7]. The special case of lattice sphere packings in three dimensions was already solved by Gauß 1831 [3] [6].

We consider the less spectacular question in the plane. In 1773 Lagrange solved the problem of the densest circle lattice packing which is indeed both an elementary and a natural question. Children imagine the right solution, and even bees know it: The densest configuration is hexagonal.

An elementary proof uses Heron´s formula for the triangle area and gives a lower bound for the area of a lattice fundamental parallelogram. An analytic-geometrical proof uses the vector product instead of Heron´s formula. Both proofs are presented and illustrated in [6]. Fukshansky describes another interesting approach using basic linear algebra und the method of successive minima [8].

For a general circle packing which is not necessarily a circle lattice packing we have the same solution. This proof is assigned to Thue 1910 and was reworked several times [9]. Zong gave a book proof in [3] using the concept of Voronoi cells and Jensen´s inequality for convex functions (also known as Höldersche Ungleichung). A more elementary approach is described in [10]: The empty-circle-property of Voronoi vertices leads to a tessellation into triangles which are estimated using trigonometric estimations.

We will examine the Voronoi cells in a saturated sphere packing which have some more useful properties compared to sole Voronoi diagrams [11] [12] [13]. So we can consider the local configuration situation around a Voronoi vertex as complementable to a lattice fundamental parallelogram configuration. For that purpose we introduce the concept of a L-triangle around a Voronoi vertex. In this way the general situation can be reduced to the lattice configuration. Thus we build a bridge from Lagrange's solution to Thue´s problem and we can give a new and elementary proof of Thue´s theorem.

# 2 Definitions and Basics

We base our definitions on [3] [14] [15] [16]. $\mathbb{R}^2$ is the two-dimensional Euclidean plane with standard norm $|.|$ . $K = \{x \in \mathbb{R}^2 : |x| \leq 1\}$ is the unit circle with area $v(K) = \pi$. Let $X$ be a set of discrete points in $\mathbb{R}^2$ with $|y - x| \geq 2$ for every pair of distinct points $x, y \in X$. Then $X$ is countable: $X = \{x_1, x_2, x_3, \dots\}$

We call (the Minkowski sum) $X + K = \{x + K; x \in X\}$ a *circle packing*.
If $X$ is a lattice

$$L = \left\{ \sum_1^2 z_i a_i \; ; \; (a_i) \text{ linearly independent}, z_i \in \mathbb{Z} \right\},$$

then we call $L + K$ a *lattice circle packing*.
In geometrical terms, $\det(L)$ is the area of the *fundamental parallelogram*

$$FP = \left\{ \sum_1^2 \lambda_i a_i \; ; \; \lambda_i \in [0; 1] \right\}$$

of $L$. A fundamental parallelogram contains exactly one circle. Thus the *circle packing density* of a given lattice circle packing can be defined as:

$$\delta(L + K) = \frac{v(K)}{det(L)}$$

Furthermore, the copies $x + FP$ ($x \in L$) yield a plane tessellation.

The *theorem of Lagrange* states that among all circle lattice packings the hexagonal circle lattice packing has the maximal packing density. An elementary proof is given in [6]. The definition of a non-lattice circle packing density is more complex [3].

For a given circle packing $X + K$ and for $x_i \in X$ we have the

**Definition 2.1.** Let $X + K$ be a circle packing configuration and let $x_i \in X$ be a configuration point. Then we call

$$V(x_i) := \left\{ x;\ |x - x_i| = \min_j |x - x_j| \right\}$$

the *Voronoi cell* of $x_i$ and $x_i$ the *center of the Voronoi cell*.

From convex geometry we know that $\bigcup_i V(x_i)$ yields a tessellation of the plane (boundary effects are not relevant) as well. Thus we have

**Definition 2.2.** Let $X + K$ be a circle packing configuration. Then we call $\bigcup_i V(x_i)$ the *Voronoi tessellation* of $X + K$.

We notice that a Voronoi tessellation of $X + K$ is unique.

The *theorem of Thue* says that a circle packing (not only a lattice circle packing) with maximal packing density must be a hexagonal packing. For a proof, the packing density of a circle packing is based on the local packing densities related to the Voronoi cells by

$$\delta(V(x_i)) = \frac{v(K)}{v(V(x_i))}$$

and it can be shown that the local inequality

$$\delta(V(x_i)) \leqq \frac{\pi}{2\sqrt{3}}$$

holds for every Voronoi cell with equality only for a regular hexagon.

A nice and compact proof of Thue´s theorem using Jensen´s inequality for convex functions is given in [3], another elementary proof using the empty-circle property of a Voronoi tessellation and trigonometric estimations can be found in [10].

Since we search for a supremum or maximum packing density, we can fill up free places in a given configuration with circles. Thus we presuppose saturated configurations in the following sense: We denote a circle packing configuration as *saturated* if it is not a strict subset of another circle packing configuration [1] [3] [10].

## 3    Basic properties of a Voronoi cell

We give a short overview of the main ideas taken from convex geometry, Voronoi diagrams and space tessellation [17] [4] [14] [13] [11] [12].

We consider a saturated infinite circle packing $X + K$. A first plane structuration of the configuration $X$ results from the so-called Voronoi diagram [11] [12] [13] or Voronoi tessellation. In this view each circle packing center $x \in X$ is assigned to a zone of influence beyond the sphere $x + K$ which has clear borders towards its nearest neighbors. For example, the Voronoi cell of a hexagonal (resp. quadratic) lattice point is a regular hexagon (resp. square) (cf. fig.1).

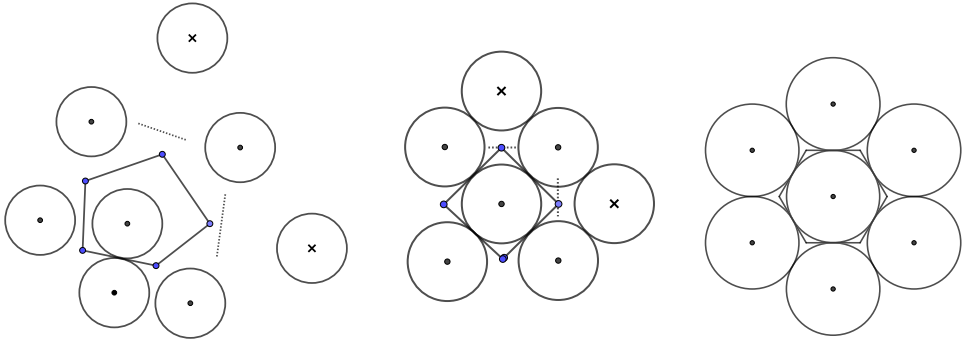

*Figure 1   Voronoi cells of a common (quadratic, hexagonal) configuration*

We consider a Voronoi cell in detail: A Voronoi cell is a polygon, its bounding edges are the perpendicular bisectors between two neighbored configuration points. Since the configuration is saturated, the volume of each Voronoi cell is finite and thus each edge is finite. Each Voronoi cell is convex, each vertex of a Voronoi cell is convex as well [4].

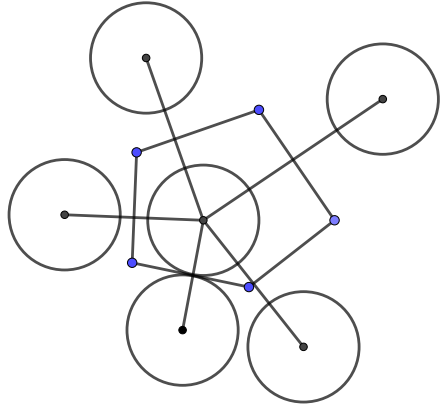

*Figure 2  Voronoi cell with perpendicular bisectors*

In particular we have

**Lemma 3.1**  Let $X + K$ be a saturated circle packing with its Voronoi tessellation and let $V(x_i)$ be a Voronoi cell with center $x_i$.
Then the nearest configuration point $x_j \in X$ generates a Voronoi edge $e$ of $V(x_i)$ and the line segment $e$ contains the intersection of the line $x_i x_j$ with the line through $e$.

**Proof.** For an elementary proof see [13, p. 59f.]. ∎

We notice that the "intersection property" of Lemma 3.1 does not hold for every Voronoi edge (cf. fig. 6). Further Voronoi cell properties are described in [13].

## 4    The situation at a Voronoi cell vertex

We sketch the so-called empty-circle-property: We consider an arbitrary center point $y \in \mathbb{R}^2$ and a circle centered in $y$ with its radius growing. This growing circle will reach for a first time and with a specific radius a point $x \in X$ of the configuration $X + K$. At this moment we have to distinguish three cases.

**Lemma 4.1.** Let $X + K$ be a saturated circle packing with its Voronoi tessellation. Let $y$ be an arbitrary point in $\mathbb{R}^2$ and $k(y)$ be the circle with center $y$ which has no point of $X$ in its interior $int(k(y))$ and at least one point of $X$ on its border $bd(k(y))$.
Case 1.  Only one point $x \in X$ is located on $bd(k(y))$. Then the circle center $y$ is located in the interior of the Voronoi cell: $y \in int(V(x_1))$
Case 2.  Two points are located on $bd(k(y))$. Then the circle center $y$ is located on the edge between two neighbored Voronoi cells.
Case 3.  Three or more points are located on $bd(k(y))$. Then the circle center $y$ is located on the vertex between three or more neighbored Voronoi cells.
The reverse is also valid in every case.

**Proof.** An elementary proof is sketched in [11, p. 3ff. und 25ff.]. ∎

Now we change the perspective and examine a circle $k(y)$ with its center $y$ located on a vertex of a Voronoi cell. Then the following Lemma holds.

**Lemma 4.2 (Empty-circle-property).** Let $X + K$ be a saturated circle packing with its Voronoi tessellation. Let $y$ be a vertex of a Voronoi cell.
Then the circle $k(y)$ through $x \in X$ contains three or more points $x_i \in X$ which span a polygon. $k(y)$ is the circumcircle of this polygon.
The interior of $k(y)$ contains no configuration point $x \in X$.
For the diameter $\text{diam}(k(y))$ holds: $\text{diam}(k(y)) < 4$

**Proof.** Applying Lemma 4.1 we get the proof for the first and the second item [11]. We can see the third item as follows: If the diameter of $k(y)$ was 4 or greater, we could place another circle at the center of the circumcircle which contradicted the assumption of a saturated circle packing configuration [10]. ∎

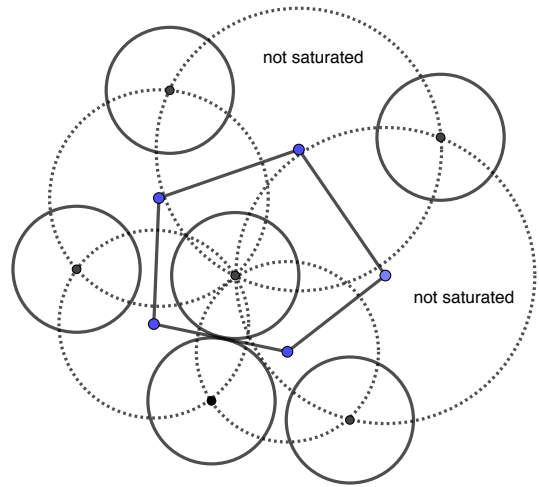

*Figure 3   Empty-circle-property*

## 5   The Voronoi vertex in a saturated circle packing

What happens if we remove a circle from a given circle packing? (cf. title fig.) We see that the Voronoi vertex expands into the gap. In particular we have the

**Lemma 5.1.** Let $X + K$ be a saturated sphere packing configuration and let $V(x_i)$ be a Voronoi cell with vertex $y_0$.
Then the distance of the Voronoi vertex from the Voronoi center is
$$|y_0 - x_i| < 2$$
and the angle between the Voronoi edges at the Voronoi vertex $\alpha$ is
$$\alpha > \frac{\pi}{3}.$$

**Proof.** Let $X + K$ be a saturated sphere packing configuration and let $V(x_i)$ be a Voronoi cell with vertex $y_0$.
We consider the limit case and assume $|y_0 - x_i| = 2$. This would be the radius of the circumcircle with centre $x_i$ according to the empty-circle-property, in contradiction to $\mathrm{diam}(k(x_i)) < 4$.
Determining the angle $\alpha$ we consider a limit case again. For the smallest possible angle we assume that every edge touches the Voronoi center circle (cf. fig. 4). And for $|y_0 - x_i| = 2$ we have by elementary trigonometry
$$\sin\left(\frac{\alpha}{2}\right) = \frac{1}{2}$$
and thus $\alpha = \frac{\pi}{3}$. From that follows the claim. ∎

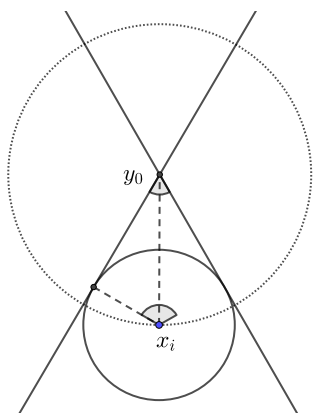

*Figure 4   Voronoi edge of a saturated circle packing*

We remain in the Voronoi vertex perspective. The question is: Can two circles around a vertex (besides the Voronoi center and the Voronoi vertex) define a fundamental parallelotope of a lattice circle packing?

## 6  Regular and degenerate Voronoi vertices

In the literature we find the usual distinction between regular and degenerated Voronoi vertices. A Voronoi vertex is called *regular* (resp. *degenerated*) if the circumcenter contains three (resp. four of more) configuration points. We give a heuristic overview of the variety of Voronoi vertices with the surrounding configuration points and polygons.

On a first glance we could assume that the regular case is simple. This is wrong, as the following examples will demonstrate: We begin with a regular vertex and intersects of the configuration points´ lines with the Voronoi cells´ edges. This case is simple, indeed. The circumcenter lies in the interior of the triangle built up by the configuration points (cf. fig. 5).

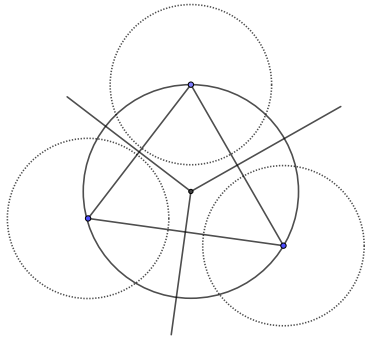

*Figure 5   Voronoi vertex - regular*

If the configuration points form a square, we have a degenerated Voronoi vertex (cf. fig. 6). If (in a first step) we shift one point away, the degenerated vertex will disintegrate into two regular vertices with a short connecting edge. The configuration points´ lines still intersect the Voronoi edges. If (in a second step) we shift two points towards the remaining quadratic Voronoi center the "vertical" configuration points´ line will not intersect the "short" Voronoi edge any more.

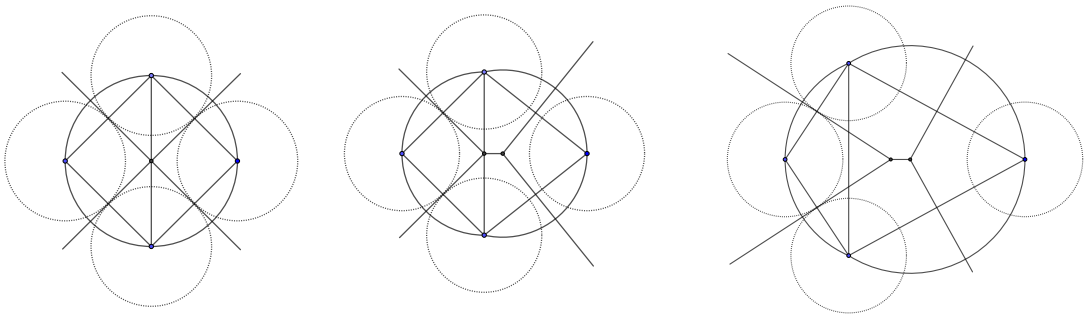

*Figure 6   Voronoi vertex – degenerated and resolved in two steps*

We can see the effect of symmetry-breaking in another way: Whereas the quadratic configuration can also be divided horizontally into two triangles, this possibility vanishes due to the broken symmetry.

We recognize: The generated "small" Voronoi edge does not cause serious problems. The quadrilateral grows with the "small" Voronoi edge. The "right" partial triangle grows as well. But the "left" partial triangle gets smaller. This phenomenon must be observed with respect of the packing density. The crucial question is: Can the configuration points´ line shifted so far to the left, that it will intersect the left circle? (cf. fig. 6)

In the literature the described triangulation phenomenon generating a "small" Voronoi edge is known as Pitteway triangulation resp. non-Pitteway triangulation. The Pitteway triangulation theorem describes the connection between a Voronoi diagram and a so-called Delone triangulation [13, p. 70ff.]. We notice that the non-Pitteway phenomenon gets relevant for Voronoi angles smaller than $\frac{\pi}{2}$. We note the Pitteway-property as

**Definition 6.1.** Let $X + K$ be a saturated sphere packing configuration with Voronoi cells $V(x_i)$. Let the Voronoi edge $e$ be generated by $x_1$ and $x_2$.
We call a Voronoi edge $e$ a Pitteway-edge if the line $x_1 x_2$ intersects $e$, otherwise a non-Pitteway edge.

Focussing on the problem of minimizing the polygon area inside a Voronoi vertex circumcircle we can concentrate on triangle situations like in fig. 5 or fig. 6 (modified).

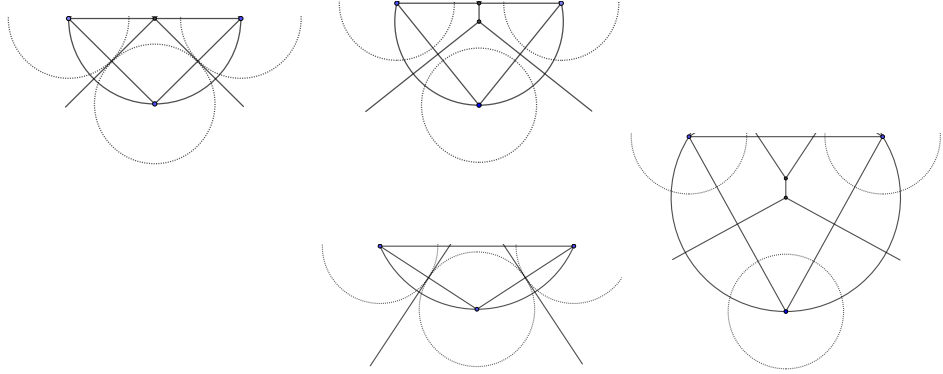

*Figure 7   Voronoi vertices – rearranged (cf. fig. 6)*

The sketched perspective will lead us to the concept of an L-triangle defined by the Voronoi center, the Voronoi vertex and the "neighbored" configuration points.

## 7   The Voronoi vertex seen from the Voronoi center

A Voronoi vertex defines a triangle in the following way: The Voronoi center is the triangle top, the neighbored circle centers form the triangle base. In general this construction is known as Delone-triangulation [12] [4] [13]. In the special case of a circle packing context with additional properties we introduce

**Definition 7.1.** Let $X + K$ be a saturated circle packing configuration and let $V(x_i)$ be a *Voronoi cell* with center $x_i$ ($x_i \in X$).
Let $y_0$ be a vertex of $V(x_i)$ and $x_{i_1}, x_{i_2} \in X$ be the symmetric configuration points of $x_i$ respecting the vertex defining edges.
Then we call the triangle defined by $x_i, x_{i_1}, x_{i_2}$ a *L-triangle* $L(y_0, x_i, x_{i_1}, x_{i_2})$.

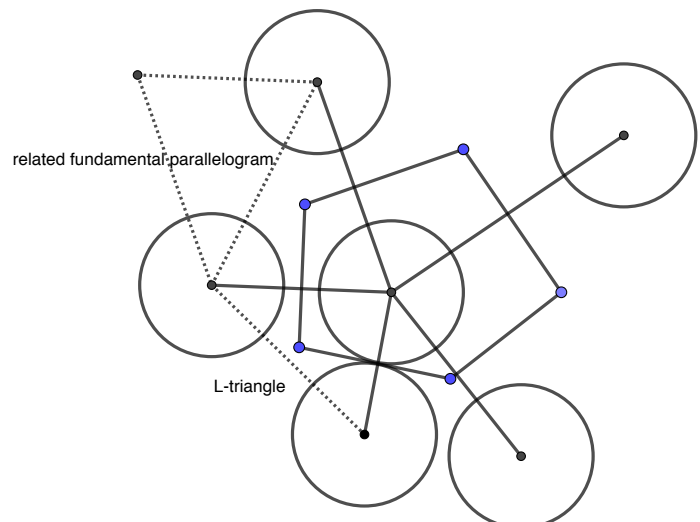

*Figure 8   L-triangle and related fundamental parallelogram*

For a circle packing we can show that the line segment $x_{i_1}x_{i_2}$ does not intersect the basic circle $x_i + K$ (cf. fig. 9).

**Theorem 7.2.** Let $X + K$ be a saturated circle packing configuration and $L(y_0, x_i, x_{i_1}, x_{i_2})$ be a L-triangle.
Then the entire circle sector of $x_i$ lies inside the L-triangle.

**Proof.** Let $L(y_0, x_i, x_{i_1}, x_{i_2})$ be a L-triangle.
Since we have a saturated packing, the diameter of the circumcircle around the vertex $y_0$ is smaller than 4 (Lemma 4.2).
We consider the limit case with circumcircle radius equal to 2 and two touching configuration circles (cf. fig. 9). According to elementary geometry we have two equilateral triangles $y_0 x_i x_{i_1}$ and $y_0 x_i x_{i_2}$. Thus we have a height of $x_i$ referred to $x_{i_1}x_{i_2}$ larger than 1. Thus the circle sector is not intersected by $x_{i_1}x_{i_2}$, which is the claim. ∎

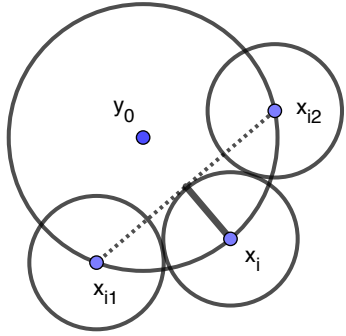

*Figure 9  L-triangle – circle sector*

**Theorem 7.3.** Let $X + K$ be a saturated circle packing configuration and $L(y_0, x_i, x_{i_1}, x_{i_2})$ be a L-triangle.
Then the L-triangle can be completed to a fundamental parallelogram of a lattice circle packing.

**Proof.** We consider an L-triangle $L(y_0, x_i, x_{i_1}, x_{i_2})$.
Let $x_i'$ be the reflected point to $x_i$ at the middle of the line segment $x_{i_1}x_{i_2}$. (In terms of analytic geometry $x_i'$ can be described as $\boldsymbol{x_i'} = \boldsymbol{x_i} + \boldsymbol{x_i x_{i_1}} + \boldsymbol{x_i x_{i_2}}$.)
Then by symmetry the circle sector around $x_i'$ is fully contained inside the parallelogram $x_i x_{i_1} x_{i_2} x_i'$. From that we can conclude [6, p. 47ff.] that a fundamental parallelogram of a circle lattice packing exists with basis $\boldsymbol{x_i x_{i_1}}$ and $\boldsymbol{x_i x_{i_2}}$. ∎

This theorem leads us to the

**Definition 7.4.** Let $X + K$ be a saturated circle packing configuration and $L(y_0, x_i, x_{i_1}, x_{i_2})$ be an L-triangle.
Then we call the parallelotope spanned by $\boldsymbol{x_i x_{i_1}}$ and $\boldsymbol{x_i x_{i_2}}$ the *related fundamental parallelogram*.

## 8   L-triangle and related fundamental parallelotope

**Theorem 8.1.** Let $X + K$ be a saturated circle packing configuration. Let $L(y_0, x_i, x_{i_1}, x_{i_2})$ be a L-triangle with its related fundamental parallelotope.
Then the area of the L-triangle will be minimal iff the area of the related fundamental parallelotope is minimal.

**Proof.** The area of the L-triangle is $\frac{1}{2}det(\boldsymbol{x_{i_1}x_{i_2}}, \boldsymbol{x_{i_1}x_i})$, the area of the related fundamental parallelotope is $det(\boldsymbol{x_{i_1}x_{i_2}}, \boldsymbol{x_{i_1}x_i})$. Since both terms are positive, the assertion is proved. ∎

## 9 Proof of Thue´s theorem

**Theorem 9.1.** Let $X + K$ be a saturated circle packing configuration.
Then the packing density will be minimal iff the circle packing configuration is a hexagonal lattice circle packing.

**Proof.** Let $X + K$ be a saturated circle packing configuration and $\bigcup_i V(x_i)$ the Voronoi tessellation of $X + K$.

Let $(y_j)$ be the set of all Voronoi vertices with circumcircles according to Lemma 4.2. We consider the polygons with circumcenter $y_j$ and configuration points $x_i, x_{i_1}, x_{i_2}, \ldots$ situated on the circumcircle.

Then we distinguish the following cases.

Case 1. The Voronoi vertex $y_j$ is the common vertex of three Voronoi cells.

We have three Voronoi centers, let them be $x_i, x_{i_1}, x_{i_2}$. They form a triangle.
According to theorem 7.3, the L-triangle $L(y_j, x_i, x_{i_1}, x_{i_2})$ can be expanded to a fundamental parallelotope, and according to theorem 8.1, the area of the L-triangle $L(y_j, x_i, x_{i_1}, x_{i_2})$ will be minimal iff the area of the related fundamental parallelotope is minimal.

Case 2. The Voronoi vertex $y_j$ is the common vertex of four Voronoi cells.

We have four Voronoi centers, let them be $x_i, x_{i_1}, x_{i_2}, x_{i_3}$. They form a quadrilateral that can be dissected into two triangles.
Let the first triangle be generated by $x_i, x_{i_1}, x_{i_2}$. Then we have an L-triangle $L(y_j, x_i, x_{i_1}, x_{i_2})$, according to the above.
Let the second triangle be generated by $x_{i_3}, x_{i_1}, x_{i_2}$. Then we have an L-triangle with the same Voronoi vertex $y_j$ but a different Voronoi center $x_{i_3}$
$L(y_j, x_{i_3}, x_{i_1}, x_{i_2})$, according to the above again.

Case 3. The Voronoi vertex $y_j$ is the common vertex of five or more Voronoi cells.

This case can be treated in the same way as case 2.

In every case we have a tessellation of the polygon with circumcenter $y_j$ and configuration points $x_i, x_{i_1}, x_{i_2}, \ldots$ into one or more L-triangles.

Since the polygons with circumcenters $y_j$ yield a tessellation of the plane, we have a further tessellation of the plane by L-triangles.

And since the L-triangles will have minimal area iff the related fundamental parallelotopes have minimal area, we obtain the proof applying Lagrange´s theorem. ■

## 10 Interpretation

The convexity of Voronoi cells plays an essential role in the solution of the Thue problem.

In the proof of Zong we can detect it in Jensen´s inequality [3]. In the given proof we find it in the convexity of a Voronoi edge and in the convexity of its circumcircle with diameter smaller than 4. Compared to Zong´s book proof, the presented proof is more elementary in the sense that analysis tools are not necessary since the Lagrange theorem can be shown only with quadratic forms [6].

As we learned from the proof, the minimum of the L-triangle area occurs for an equilateral triangle. This observation corresponds to the fundamental parallelogram of a hexagonal lattice as well as to the sexangular Voronoi cell of a hexagonal lattice.

And what is astonishing, the minimum of the L-triangle also occurs as limit value of an L-triangle with 120°-angle. We recognize that this L-triangle arises when a fundamental parallelogram of a hexagonal lattice is divided in an alternative manner. Solely because the diameter of the circumcircle cannot equal 4, this phenomenon loses its relevance for the presented proof.

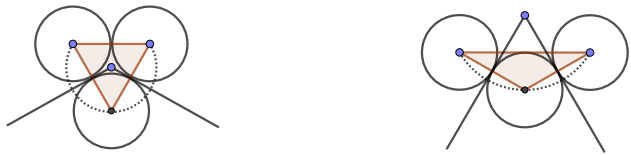

*Figure 10   L-triangle - minimum and infimum*

**Acknowledgements**
I owe many thanks to V. Ulm for his support, A. Beutelspacher for his encouragement, Eric Müller for useful mathematical hints and Th. Zeitlhöfler for proof reading.